\documentclass[preprint,12pt]{elsarticle}
\usepackage{amsfonts, amssymb, amsmath, mathrsfs}
\usepackage[top=1in,bottom=1in,left=1.25in,right=1.25in]{geometry}
\journal{{\rm } }

\begin{document}
\newtheorem{lemma}{Lemma}[section]
\newtheorem{proposition}{Proposition}[section]
\newtheorem{theorem}{Theorem}[section]
\newtheorem{corollary}{Corollary}[section]
\newtheorem{example}{Example}[section]
\newtheorem{definition}{Definition}[section]
\newtheorem{remark}{Remark}[section]
\newtheorem{property}{Property}[section]
 \makeatletter
    \newcommand{\rmnum}[1]{\romannumeral #1}
    \newcommand{\Rmnum}[1]{\expandafter\@slowromancap\romannumeral #1@}
    \makeatother

\begin{frontmatter}
\title{Comparison of Two Types of Separation Axioms in Soft Topological Spaces}

 \author[label1]{Guan Xuechong\corref{cor1}}
 \ead{guanxc@foxmail.com}
 \cortext[cor1]{School of Mathematics and Statistics, Jiangsu Normal University, Xuzhou, 221116, China. }
\address[label1]{School of Mathematics and Statistics, Jiangsu Normal University, Xuzhou, 221116, China}

\begin{abstract}
In the study of soft topological spaces, two types of separation axioms have given if soft points and single
point soft sets have been taken as separated objects respectively.
In this paper, some examples and properties are given to explore differences between
separation axioms given in \cite{topological,Separation}.
\end{abstract}


\end{frontmatter}

\section{Introduction }\label{intro}
In order to solve the problem of data uncertainty, mathematicians have put forward research schemes such as fuzzy set theory,
interval mathematics theory and rough set theory. However, these theories also have their own defects.
In the year 1999, Molodtsov initiated soft set theory \cite{first}
which is a mathematical tool for parameterizing research objects of sets.
At present this theory has been widely combined with the theories of fuzzy set \cite{fuzzy}-\cite{interval},
algebra \cite{category}-\cite{rough}, topology \cite{topology}-\cite{zhang}, and decision-making
\cite{reduction}-\cite{making},
and rich results have been obtained.

In the research of soft topological spaces,
soft separation axioms which are draw lessons from general topology theory are given.
The most remarkable thing among them is the difference of separated objects.
In \cite{topological} the separated objects are single point soft sets, such as $\widetilde{\{x\}}$.
Soft points denoted by $e_F$ and $P_e^x$ respectively in \cite{Sabir,Separation} are also regarded
as the separated objects in separation axioms.
Here we adopt the definition of soft points given in \cite{Separation}.
The main purpose of this paper is to compare the relationship between two types of separation axioms in soft
topological spaces.

The rest of this paper is organized as follows.
The second section briefly reviews some basic notions on soft theory.
The third section lists two types of separation axioms in soft topological spaces,
and then compares the differences among them by examples.
Some properties on the relationship of separation axioms are also shown.

\section{Preliminaries}
\label{sec:2}

In this section, first, we present some basic definitions and results of soft set theory.
Let $X$ be an initial universe set and $E$ be a set of parameters which usually are initial attributes,
characteristics, or properties of objects in $X$.

\begin{definition}\label{definition:1} {\rm (\cite{first})}
{\rm
Let $X$ be an initial universe set and $E$ be a set of parameters. Let ${\cal P}(X)$ denote the power set of $X$ and
$A\subseteq E$.
A pair $(F,A)$ is called a soft set over $X$, where $F$ is a mapping given by $F : A\rightarrow {\cal P}(X)$.}
\end{definition}

In other words, a soft set over $X$ is a parameterized family of subsets of the universe $X$. For $\varepsilon\in A$,
$F(\varepsilon)$ may be considered as the set of $\varepsilon$-approximate elements of the soft set $(F,A)$.
The family of all soft sets over X with a parameters set $E$ is denoted by $S\!S_E(X)$.

\begin{definition}\label{definition:2}{\rm (\cite{topological})}
{\rm Let $Y$ be a non-empty subset of $X$, then $\widetilde{Y}$ denotes the soft set $(Y, E)$ over $X$
 for which $Y(e)= Y$, for all $e\in E$.
If $Y$ is a single point set, such as $Y=\{x\}$ and $x\in X$, then
$\widetilde{Y}$ is called a single point soft set.

}
\end{definition}

\begin{definition}\label{definition:3}{\rm (\cite{rough})}
{\rm For two soft sets $(F,A)$ and $(G,B)$ over a common universe $X$,
$(F, A)$ is a soft subset of $(G,B)$, denoted by $(F, A)\subseteq (G, B)$, if

(i) $A \subseteq B$, and

(ii) $\forall \varepsilon\in A$, $F(\varepsilon)$ is a subset of $G(\varepsilon)$.

If $(F, A)\subseteq (G, B)$ and $(G, B)\subseteq (F, A)$, we write $(F, A)=(G, B)$
and say $(F, A)$ is equal to $(G, B)$.
}
\end{definition}

\begin{definition}\label{definition:4}{\rm (\cite{topological})}
{\rm Let $(F,E)$ be a soft sets over a common universe $X$ and $x\in X$.
We say that $x \widehat{\in} (F , E)$ read as $x$ belongs to the soft set $(F , E)$
whenever $x\in F(e)$ for all $e\in E$.
}
\end{definition}

\begin{definition}\label{definition:7}{\rm (\cite{soft})}
{\rm The union of two soft sets $(F, A)$ and $(G, B)$
over a common universe $X$ is the soft set $(H, C)$, where $C = A \cup B$, and $\forall e \in C$,
$$H(e)=\left\{
          \begin {array}{ll}
          F(e),&{\mbox{if}} \ \ e\in A-B;\\
          G(e),&{\mbox{if}} \ \ e\in B-A;\\
          F(e)\cup G(e),&{\mbox{if}} \ \ e\in A\cap B.
         \end{array}
        \right.$$
We write $(F,A)\cup (G,B)=(H,C)$.}
\end{definition}

For the union of a family of soft sets $\{(F_i, A_i): i\in I\}$,
$$ \bigcup\limits_{i\in I} (F_i, A_i)=(F, \bigcup\limits_{i\in I} A_i),$$
where $F(e) = \bigcup \{ {F_i}(e):i \in {I_e}\}$ and ${I_e} = \{ i \in I: e \in {A_i}\}$
(see \cite{semirings,continuous}).

\begin{definition}\label{definition:11}{\rm (\cite{operation})}
{\rm The intersection of two soft sets $(F, A)$ and $(G, B)$ over a common universe $X$ such that $A\cap B\neq \emptyset$
is the soft set $(H, C)$, where $C = A \cap B$, and $\forall e \in C$,
$H(e)= F(e)\cap G(e).$ We write $(F,A)\cap (G,B)=(H,C)$.}
\end{definition}

\begin{definition}\label{definition:12}{\rm (\cite{topological})}
{\rm Let ${\cal T}$ be the collection of soft sets over $X$, then ${\cal T}$ is said to be a soft topology on $X$ if\\
(1) $\widetilde{X},\widetilde{\emptyset}$ belong to ${\cal T}$\\
(2) the union of any number of soft sets in ${\cal T}$ belongs to ${\cal T}$\\
(3) the intersection of any two soft sets in ${\cal T}$ belongs to ${\cal T}$.\\
The triplet $(X, {\cal T}, E)$ is called a soft topological space over $X$.
The members of ${\cal T}$ are said to be soft open sets in $X$.
}
\end{definition}

The definition of  soft topological spaces means the family of soft sets in a soft topology is closed
with respect to the operations of union and finite intersection.

\begin{example}\label{example:1}
{\rm Let $X=\{x, y\}$, $E=\{e_1,e_2\}$ and
${\cal T}=\{ \widetilde \emptyset ,\widetilde X,({F_1},E),({F_2},E)\}$, where
$${F_1}(e_1) = \left\{ x \right\},F_1(e_2) =
\left\{ y \right\};{\kern 1pt} {\kern 1pt} F_2(e_1)
= \left\{ y \right\},F_2(e_2) = \left\{ x \right\}$$
Clearly, $(X, {\cal T}, E)$ is a soft topological space over $X$.
}
\end{example}

In the following we give two examples of non-soft topological spaces.

\begin{example}\label{example:2}
{\rm Let $X=\{x, y\}$, $E=\{e_1,e_2\}$ and
${\cal T}=\{\widetilde \emptyset ,\widetilde X,(F_1,E ),(F_2,E ),(F_3,E),\\
(F_4, E), (F_5, E), (F_6, E), (F_7, E), (F_8,E),
(F_9, E), (F_{10}, E), (F_{11}, E)\}$,
where\\
${F_1}({e_1}) = \{x \},{F_1}({e_2})=\{y\}$,\\
${F_2}({e_1}) = \{ y \},{F_2}({e_2}) = \emptyset $,\\
${F_3}({e_1}) = \emptyset ,{F_3}({e_2}) = \{ x \}$,\\
${F_4}({e_1}) = \{ y \},{F_4}({e_2}) =\{y \}$,\\
${F_5}({e_1}) = \{ x \},{F_5}({e_2}) = \emptyset$,\\
${F_6}({e_1}) = X,{F_6}({e_2}) = \{ y \}$,\\
${F_7}({e_1}) = \{ x \},{F_7}({e_2}) = X$,\\
${F_8}({e_1}) = \{ y \},{F_8}({e_2}) = \{ x \}$,\\
${F_9}({e_1}) = X,{F_9}({e_2}) = \emptyset$,\\
${F_{10}}({e_1}) = \{ y \},{F_{10}}({e_2}) = X$,\\
${F_{11}}({e_1}) = \{ x \},{F_{11}}({e_2}) = \{x \}.$

The triple $(X, {\cal T}, E)$ is given in \cite{Separation}. But it is not a soft topological space over $X$, since
$({F_8},E) \cup ({F_9},E)\not\in {\cal T},({F_7},E) \cap ({F_{10}},E)\not\in {\cal T}, ({F_4},E) \cap ({F_7},E)\not\in {\cal T}$.
Assume that three mappings $F_{12}, F_{13},F_{14}$ are defined as:\\
${F_{12}}({e_1}) = X,{F_{12}}({e_2}) =\{ x \} $,\\
${F_{13}}({e_1}) = \emptyset,{F_{13}}({e_2}) = X$,\\
${F_{14}}({e_1}) = \emptyset,{F_{14}}({e_2}) = \{y\}.$\\
Let ${\cal T}^* ={\cal T}\cup\{(F_{12},E),(F_{13},E),(F_{14},E)\}= S\!S_E(X)$.
Then $(X, {\cal T}^*, E)$ is a soft topological space.
}
\end{example}

\begin{example}\label{example:3}
{\rm Let $X = \{ {x_1},{x_2},{x_3},{x_4}\}$, $E = \left\{ {{e_1},{e_2}} \right\}$
and
${\cal T} = \{ \widetilde \emptyset ,\widetilde X,\left( {{F_1},E} \right),\left( {{F_2},E} \right),\\
\left( {{F_3},E} \right),
\left( {{F_4},E} \right),\left( {{F_5},E} \right),\left( {{F_6},E} \right),\left( {{F_7},E} \right),
\left( {{F_8},E} \right),\left( {{F_9},E} \right),\left( {{F_{10}},E} \right),\left( {{F_{11}},E} \right),
\left( {{F_{12}},E} \right),\\\left( {{F_1}_3,E} \right),
\left( {{F_1}_4,E} \right),\left( {{F_{15}},E} \right),\left( {{F_{16}},E} \right),\left( {{F_{17}},E} \right),
\left( {{F_{18}},E} \right),\left( {{F_{19}},E} \right),
\left( {{F_{20}},E} \right),\left( {{F_{21}},E} \right),\\
\left( {{F_{22}},E} \right)\}$, where\\
${F_1}({e_1}) = \{ {x_1},{x_2},{x_4}\} ,{F_1}({e_2}) = \{ {x_1},{x_2},{x_3}\}$,\\
${F_2}({e_1}) = \{ {x_1},{x_3},{x_4}\} ,{F_2}({e_2}) = \{ {x_1},{x_2},{x_3}\}$,\\
${F_3}({e_1}) = \{ {x_1},{x_4}\} ,{F_3}({e_2}) = \{ {x_1},{x_2},{x_3}\}$,\\
${F_4}({e_1}) = \{ {x_2},{x_3}\} ,{F_4}({e_2}) = \{ {x_1},{x_2},{x_3}\}$,\\
${F_5}({e_1}) = \left\{ {{x_2}} \right\},{F_5}({e_2}) = \{ {x_1},{x_2},{x_3}\}$,\\
${F_6}({e_1}) = \left\{ {{x_3}} \right\},{F_6}({e_2}) = \{ {x_1},{x_2},{x_3}\}$,\\
${F_7}({e_1}) = \emptyset ,{F_7}({e_2}) = \{ {x_1},{x_2},{x_3}\}$,\\
${F_8}({e_1}) = \{ {x_2},{x_3},{x_4}\} ,{F_8}({e_2}) = \{ {x_1},{x_4}\}$,\\
${F_9}({e_1}) = \left\{ {{x_1}} \right\},{F_9}({e_2}) = \{ {x_2},{x_3}\}$,\\
${F_{10}}({e_1}) = \{ {x_2},{x_4}\} ,{F_{10}}({e_2}) = \left\{ {{x_1}} \right\}$,\\
${F_{11}}({e_1}) = \{ {x_3},{x_4}\} ,{F_{11}}({e_2}) = \left\{ {{x_1}} \right\}$,\\
${F_{12}}({e_1}) = \left\{ {{x_4}} \right\},{F_{12}}({e_2}) = \left\{ {{x_1}} \right\}$,\\
${F_{13}}({e_1}) = \{ {x_2},{x_3}\} ,{F_{13}}({e_2}) = \left\{ {{x_1}} \right\}$,\\
${F_{14}}({e_1}) = \left\{ {{x_2}} \right\},{F_{14}}({e_2}) = \left\{ {{x_1}} \right\}$,\\
${F_{15}}({e_1}) = \left\{ {{x_3}} \right\},{F_{15}}({e_2}) = \left\{ {{x_1}} \right\}$,\\
${F_{16}}({e_1}) = \emptyset ,{F_{16}}({e_2}) = \left\{ {{x_1}} \right\}$,\\
${F_{17}}({e_1}) = \{ {x_2},{x_3},{x_4}\} ,{F_{17}}({e_2}) = X$,\\
${F_{18}}({e_1}) = \emptyset ,{F_{18}}({e_2}) = \{ {x_2},{x_3}\}$,\\
${F_{19}}({e_1}) = \{ {x_1},{x_2},{x_3}\} ,{F_{19}}({e_2}) = \{ {x_1},{x_2},{x_3}\}$,\\
${F_{20}}({e_1}) = \{ {x_1},{x_2}\} ,{F_{20}}({e_2}) = \{ {x_1},{x_2},{x_3}\}$,\\
${F_{21}}({e_1}) = \{ {x_1},{x_3}\} ,{F_{21}}({e_2}) = \{ {x_1},{x_2},{x_3}\}$,\\
${F_{22}}({e_1}) = \left\{ {{x_1}} \right\},{F_{22}}({e_2}) = \{ {x_1},{x_2},{x_3}\}$.

The triple $(X, {\cal T}, E)$ is given in \cite{zhang}. But it is not a soft topological space over $X$, since
$({F_3},E) \cup ({F_4},E),({F_2},E) \cup ({F_4},E), ({F_4},E) \cup ({F_8},E),
({F_4},E) \cup ({F_{12}},E),({F_8},E) \cup ({F_{12}},E),({F_{11}},E) \cup ({F_{13}},E)$ are not in ${\cal T}$.

Suppose that the mappings $F_{23}, F_{24}, F_{25}, F_{26}, F_{27}$ are defined as follows:\\
${F_{23}}({e_1}) = X,{F_{23}}({e_2}) = \{ {x_1},{x_2},{x_3}\}\\
 {F_{24}}({e_1}) = \{ {x_2},{x_3},{x_4}\} ,{F_{24}}({e_2}) = \{ {x_1}\}$,\\
${F_{25}}({e_1}) = \{ {x_2},{x_3},{x_4}\} ,{F_{25}}({e_2}) = \{ {x_1},{x_2},{x_3}\}\\
 {F_{26}}({e_1}) = \{ {x_2},{x_3},{x_4}\} ,{F_{26}}({e_2}) = X$,\\
${F_{27}}({e_1}) = \{ {x_2},{x_3},{x_4}\} ,{F_{27}}({e_2}) = \{ {x_1},{x_4}\}$.\\
Let ${{\cal T}^{*}}{\rm{ = }}{\cal T} \cup \{ ({F_{23}},E),({F_{24}},E),({F_{25}},E),({F_{26}},E),({F_{27}},E)\}$.
Then $(X, {\cal T}^*, E)$ is a soft topological space over $X$.
}
\end{example}

\begin{definition}\label{definition:5}{\rm (\cite{topological})}
{\rm
The relative complement of a soft set $(F, A)$ is denoted by $(F, A)^{'}$ and is defined by $(F, A)^{'} = (F^{'}, A)$ where
$F^{'}: A \rightarrow {\cal P}(X)$ is a mapping given by
$F^{'}(e) = X- F(e)$ for all $e\in E$.
}
\end{definition}

\begin{definition}\label{definition:6}{\rm (\cite{topological})}
{\rm
Let $(X, {\cal T}, E)$ be a soft topological space over $X$. If $(F, E)^{'}\in {\cal T}$,
then $(F, E)$ is called a soft closed set.
}
\end{definition}

\begin{definition}\label{definition:8}{\rm (\cite{Separation})}
{\rm
A soft set $(P, E)$ is called a soft point, denoted by $P_{e_0}^{x_0}$,
if there exist $e_0\in E$ and $x_0\in X$ such that $P(e_0)=\{x_0\}$ and $P(e)=\emptyset$ for all $e\in E-\{e_0\}$.
The family of all soft points in $S\!S_E(X)$ is denoted by $S\!P(S\!S_E(X))$.
}
\end{definition}

It should be noted that a soft point defined here is different from \cite{remark}.
For a soft point $P_{e_0}^{x_0}$ and a soft set $(G,E)$, we write $P_{e_0}^{x_0}\widetilde{\in} (G,E)$
if $P_{e_0}^{x_0}\subseteq (G,E)$(see \cite{remark}).
It is clearly that $P_{e_0}^{x_0}\widetilde{\in} (G,E) \Leftrightarrow x_0\in G(e_0)$.
For another symbol $\widehat{\in}$, there exists a similar relationship that
$x\widehat{\in} (F,E)\Leftrightarrow \widetilde{\{x\}}\subseteq (F,E)$.
We have $P_{e_0}^{x_0}\widetilde{\in} (G,E)$ if $x_0\widehat{\in} (G, E)$. But the converse is generally not true.
Essentially the relationships $\widehat{\in}, \widetilde{\in}$
between soft sets are all inclusion relations $\subseteq$.

By introducing soft points and single point soft sets,
then the soft separation axioms appear in soft topological spaces naturally.

\section{Separation axioms in Soft topological spaces}
\label{sec:3}

Separation axioms are important topological properties. It is of practical significance to discuss some relevant separation
axioms in soft set theory.
In a sense it indicates the information represented by soft sets is incompatible,
if the intersection of two soft sets is empty.
The separation axioms of soft topological space further explains, for any two different basic elements (such as soft points),
incompatible soft sets with more information can be found.

At present, there are different ways to take ``points" in the separation of soft topological spaces.
References \cite{topological,Sabir,Separation}
took single point soft sets or soft points generated by elements in the initial set as separation objects respectively.
In this paper in order to distinguish the definitions of soft separations in \cite{topological} and \cite{Separation},
these same names are marked differently.
The following sections are treated similarly.

\subsection{Soft $T_0-$\Rmnum{1} space and soft $T_0$-\Rmnum{2} space}

\begin{definition}\label{definition:9}{\rm (\cite{topological})}
{\rm
Let $(X, {\cal T}, E)$ be a soft topological space over $X$ and $x, y\in X$ such that $x\neq y$.
If there exist soft open sets $(F , E)$ and $(G, E)$ such that
$x\widehat{\in} (F , E)$ and $y\widehat{\not\in} (F , E)$ or
$y\widehat{\in} (G, E)$ and $x\widehat{\not\in} (G, E)$, then $(X, {\cal T}, E)$ is called a soft $T_0$-\Rmnum{1} space.
}
\end{definition}

\begin{definition}\label{definition:10}{\rm (\cite{Separation})}
{\rm
Let $(X, {\cal T}, E)$ be a soft topological space over $X$. $P_a^x$ and $P_b^y$ are two different soft points in $S\!P(S\!S_E(X))$.
If there exist soft open sets $(F , E)$ and $(G, E)$ such that
$P_a^x\widehat{\in} (F , E)$ and $P_b^y\widehat{\not\in} (F , E)$ or
$P_b^y\widehat{\in} (G, E)$ and $P_a^x\widehat{\not\in} (G, E)$, then $(X, {\cal T}, E)$ is called a soft $T_0$-\Rmnum{2} space.
}
\end{definition}

Obviously separability of soft topological spaces means objects such as $P_a^x$ and $\widetilde{\{x\}}$
can be separated by soft open sets.
For two finite sets $X$ and $E$, the cardinality of $S\!P(S\!S_E(X))$ is generally bigger than that of the set $X$.
For example, let $| X |=m$ and $| E |=n$, then $|S\!P(S\!S_E(X))|=m\times n$.
Therefore, it is a bit more complicated to verify the separation axioms $T_0$-\Rmnum{2} than the separation $T_0$-\Rmnum{1}.

A soft topological space does not necessarily have the above separation axioms. It is shown in the following examples.

\begin{example}\label{example:4}
{\rm
Let $(X, {\cal T}, E)$ be a soft topological space over $X$ as shown in Example \ref{example:1}.
For all $x\in X$, if $x\widehat\in (F,E)$, then $(F,E)=\widetilde X$. While $y\widehat\in \widetilde X$.
It indicates that $x$ and $y$ cannot be separated by soft open sets. So $(X, {\cal T}, E)$ is not a soft $T_0$-\Rmnum{1} space.

Take $P_{{e_1}}^{y},P_{{e_2}}^{x} \in S\!P(S\!S_E(X))$. We have $P_{{e_1}}^{y} \ne P_{{e_2}}^{x}$.
The soft set $(F_2,E)$ is the only one nontrivial soft open set such that $P_{{e_1}}^{y}\widetilde  \in (F_2,E)$
or $P_{{e_2}}^{x}\widetilde  \in (F_2,E)$.
Thus $(X, {\cal T}, E)$ is also not a soft $T_0$-\Rmnum{2} space.
}
\end{example}

\begin{example}\label{example:5}
{\rm
Let $X=\{x, y\}$, $E=\{e_1,e_2\}$ and
${\cal T}=\{\widetilde \emptyset ,\widetilde X,(F,E )\}$, where
$$F(e_1)=F(e_2)=\{x\}.$$
Then $(X, {\cal T}, E)$ is a soft topological space.

For $x\neq y$, we have $x\widehat\in (F,E)$ and $y\widehat{\not\in} (F,E)$.
So $(X, {\cal T}, E)$ is a soft $T_0$-\Rmnum{1} space.

Take $P_{{e_1}}^{x},P_{{e_2}}^{x} \in S\!P(S\!S_E(X))$. We have $P_{{e_1}}^{x} \ne P_{{e_2}}^{x}$.
While $(F,E)$ is the only one nontrivial soft open set such that $P_{{e_1}}^{x}\widetilde  \in (F,E)$
or $P_{{e_2}}^{x}\widetilde  \in (F,E)$.
Thus $(X, {\cal T}, E)$ is not a soft $T_0$-\Rmnum{2} space.

}
\end{example}

\begin{remark}\label{remark:1}
{\rm
Example \ref{example:5} shows that a soft $T_0$-\Rmnum{1} space may not to be a soft $T_0$-\Rmnum{2} space.
}
\end{remark}

\begin{proposition}\label{proposition:5}
{\rm
Let $X=\{x, y\}$, $E=\{e_1,e_2\}$. If a soft topology $\cal T$ over $X$ and $E$ contains any of the following six soft sets
$(F_i, E)(i=1,2,\cdots,6)$, then $(X, {\cal T}, E)$ is a soft $T_0$-\Rmnum{1} space, where\\
${F_1}({e_1}) = \{x \},{F_1}({e_2})=\{x\}$,\\
${F_2}({e_1}) = \{ y \},{F_2}({e_2}) = \{y\}$,\\
${F_3}({e_1}) = \{x \},{F_3}({e_2}) = X $,\\
${F_4}({e_1}) = \{ y \},{F_4}({e_2}) =X $,\\
${F_5}({e_1}) = X,{F_5}({e_2}) = \{x \}$,\\
${F_6}({e_1}) = X,{F_6}({e_2}) = \{ y \}$.
}
\end{proposition}
\noindent {\bf Proof.} For all $i\in I$, we have $x\widehat{\in} (F_i, E)$ and $y\widehat{\not\in} (F_i, E)$ or
$y\widehat{\in} (F_i, E)$ and $x\widehat{\not\in} (F_i, E)$.
Then, by Definition \ref{definition:9}, we obtain the conclusion straightforwardly.
\qed

\begin{proposition}\label{proposition:1}
{\rm
Let $X=\{x, y\}$, $E=\{e_1,e_2\}$.
If $(X, {\cal T}, E)$ is a soft $T_0$-\Rmnum{2} space, then $(X, {\cal T}, E)$ is a soft $T_0$-\Rmnum{1} space.
}
\end{proposition}
\noindent {\bf Proof.}
(Proof by contradiction) Assume that $(X, {\cal T}, E)$ is not a soft $T_0$-\Rmnum{1} space.
We denote ${\cal B}=\{(F_i,E): i=1,2,\cdots,6\}$, where $F_i$ is defined as given in Proposition \ref{proposition:5}.
Then ${\cal B}\cap {\cal T}=\emptyset$.

We have $S\!P(S\!S_E(X))=\{P_{e_1}^{x}, P_{e_1}^{y}, P_{e_2}^{x}, P_{e_2}^{y}\}$.
These four soft sets form six pairs of soft points as follows:
$${\rm (\rmnum{1})} P_{e_1}^{x}, P_{e_1}^{y} \ \
{\rm (\rmnum{2})} P_{e_1}^{x}, P_{e_2}^{x}\ \
{\rm (\rmnum{3})} P_{e_1}^{x}, P_{e_2}^{y}\ \
{\rm (\rmnum{4})}  P_{e_1}^{y}, P_{e_2}^{x}\ \
{\rm (\rmnum{5})}  P_{e_1}^{y}, P_{e_2}^{y}\ \
{\rm (\rmnum{6})}  P_{e_2}^{x}, P_{e_2}^{y}.
$$
Because $(X, {\cal T}, E)$ is a soft $T_0$-\Rmnum{2} space and $P_{e_1}^{x}\neq P_{e_1}^{y}$,
the soft topology ${\cal T}$ contains at least one of the soft sets $P_{e_1}^{x}, P_{e_1}^{y},(F_7,E)$ and $(F_8,E)$, where\\
${F_7}({e_1}) = \left\{ {x} \right\},{F_7}({e_2}) = \left\{y \right\},$\\
${F_8}({e_1}) = \left\{ y \right\},{F_8}({e_2}) = \left\{ x \right\}.$\\
For the group (\rmnum{3}), we have $P_{e_1}^{x}\neq P_{e_2}^{y}$.
Then ${\cal T}$ contains at least one of the soft sets $P_{e_1}^{x}, P_{e_2}^{y}, (F_9,E)$ and $(F_{10},E)$, where\\
${F_9}({e_1}) = X,{F_9}({e_2}) = \emptyset, $\\
${F_{10}}({e_1}) = \emptyset,{F_{10}}({e_2}) = X.$
\\
For the group (\rmnum{4}), we have $P_{e_1}^{y}\neq P_{e_2}^{x}$.
Then ${\cal T}$ contains at least one of the soft sets $P_{e_1}^{y},P_{e_2}^{x},(F_9,E)$ and $(F_{10},E)$.
For the group (\rmnum{6}), we have $P_{e_2}^{x}\neq P_{e_2}^{y}$.
Then ${\cal T}$ contains at least one of the soft sets $P_{e_2}^{x},P_{e_2}^{y},(F_7,E)$ and $(F_8,E)$.

Furthermore, soft sets for separating these soft points are further determined below according to
the closure of $\cal T$ on the union operation. If the soft set $(F_9,E)$ separating the group
(\rmnum{4}) and any soft set separating the group (\rmnum{6}), or
$(F_{10},E)$ separating the group (\rmnum{4}) and any soft set
separating the group (\rmnum{1}) belong to a same topology, there
must have a soft set $(F_i,E)$ in the set $\cal B$.
For example, if $(F_9,E)$ and $P_{e_2}^{x}$ belong to ${\cal T}$, then $(F_9,E)\cup P_{e_2}^{x}=(F_5,E)\in {\cal T}$.
It is contradictory to ${\cal B}\cap {\cal T}=\emptyset$.
So $P_{e_1}^{y}$ or $P_{e_2}^{x}$ is the only soft sets separating the group (\rmnum{4}).

If $P_{e_1}^{y}\in {\cal T}$, then the soft set separating the group (\rmnum{6}) is $P_{e_2}^{x}$ or $(F_8,E)$.
If $P_{e_2}^{x}\in {\cal T}$, then $P_{e_1}^{y}\cup P_{e_2}^{x}=(F_8,E)\in {\cal T}$.
That is to say, we always have $(F_8,E)\in {\cal T}$. But it is clearly that $(F_8,E)$ and any soft set separating
the group (\rmnum{3}) cannot exist in a same topology which is not a soft $T_0$-\Rmnum{1} topology.
Thus the group (\rmnum{3}) cannot be separated.

If $P_{e_2}^{x}\in {\cal T}$, then $P_{e_2}^{y}$ or $(F_{10},E)$ is the soft set separating the group (\rmnum{3}).
If $P_{e_2}^{y}\in {\cal T}$, then $P_{e_2}^{y}\cup P_{e_2}^{x}=(F_{10},E)\in {\cal T}$.
Thus we always have $(F_{10},E)\in {\cal T}$.
While the soft set $(F_{10},E)$ and any soft set separating
the group (\rmnum{1}) cannot exist in a same topology which is not a soft $T_0$-\Rmnum{1} topology.
Then the group (\rmnum{1}) cannot be separated.

So the assumption is not true in any case. That is to say, $(X, {\cal T}, E)$ is a soft $T_0$-\Rmnum{1} space.
\qed


\subsection{Soft $T_1-$\Rmnum{1} space and soft $T_1$-\Rmnum{2} space}
\begin{definition}\label{definition:13}{\rm (\cite{topological})}
{\rm
Let $(X, {\cal T}, E)$ be a soft topological space over $X$ and $x, y\in X$ such that $x\neq y$.
If there exist soft open sets $(F , E)$ and $(G, E)$ such that
$x\widehat{\in} (F , E)$ and $y\widehat{\not\in} (F , E)$,
$y\widehat{\in} (G, E)$ and $x\widehat{\not\in} (G, E)$, then $(X, {\cal T}, E)$ is called a soft $T_1$-\Rmnum{1} space.
}
\end{definition}

\begin{definition}\label{definition:14}{\rm (\cite{Separation})}
{\rm
Let $(X, {\cal T}, E)$ be a soft topological space over $X$ and $P_a^x, P_b^y\in S\!P(S\!S_E(X))$ such that $P_a^x\neq P_b^y$.
If there exist soft open sets $(F , E)$ and $(G, E)$ such that
$P_a^x\widetilde{\in} (F , E)$ and $P_b^y\widetilde{\not\in} (F , E)$,
$P_b^y\widetilde{\in} (G, E)$ and $P_a^x\widetilde{\not\in} (G, E)$, then $(X, {\cal T}, E)$ is called a soft $T_1$-\Rmnum{2} space.
}
\end{definition}

Firstly, a soft topological space satisfying separation axioms $T_1$-\Rmnum{1} and $T_1$-\Rmnum{2} is constructed.

\begin{example}\label{example:6}
{\rm Let $X=\{x\}$, $E=\{e_1,e_2\}$ and
${\cal T}=\{ \widetilde \emptyset ,\widetilde X, P_{{e_1}}^x, P_{{e_2}}^x\}$.
The triple $(X, {\cal T}, E)$ is a soft topological space over $X$.
Since $X$ is a single point set, $(X, {\cal T}, E)$ is a soft $T_1$-\Rmnum{1} topological space.

Here $S\!P(S\!S_E(X))= \{P_{{e_1}}^x, P_{{e_2}}^x\}$.
We have
$P_{{e_1}}^x \widetilde\in P_{e_1}^x, P_{{e_1}}^x \widetilde  \notin P_{e_2}^x$
and $P_{{e_2}}^x \widetilde\in P_{e_2}^x, P_{{e_2}}^x \widetilde  \notin P_{{e_1}}^x$.
Thus $P_{{e_1}}^x$ and $P_{{e_2}}^x$ can be separated by themselves. Then $(X, {\cal T}, E)$ is a soft $T_1$-\Rmnum{2} space.
}
\end{example}

\begin{example}\label{example:7}
{\rm Let $X=\{x, y\}$, $E=\{e_1,e_2\}$ and
${\cal T}=\{\widetilde \emptyset ,\widetilde X,(F,E ),(G,E )\}$,
where
$${F}({e_1}) = {F}({e_2})=\{x\},\ \ G({e_1}) =G({e_2})= \{ y \}$$
For $x$ and $y$, we have $x\widehat\in (F,E)$, $y\widehat\notin (F,E)$
and $y\widehat\in (G,E)$, $x\widehat\notin (G,E)$.
Then $(X, {\cal T}, E)$ is a soft $T_1$-\Rmnum{1} topological space.

Let $(H,E)\in {\cal T}$ be a nontrivial soft set. If $P_{e_1}^ x\widetilde\in (H,E)$,
then $x\in H(e_1)$. We have $(H,E)=(F,E)$. Thus $P_{e_2}^ x\widetilde\in (H,E)$.
That is to say, if $P_{e_1}^ x\widetilde\in (H,E)$,
then $P_{e_2}^ x\widetilde\in (H,E)$.
The converse is also true.
So $(X, {\cal T}, E)$ is not a soft $T_1$-\Rmnum{2} topological space.
}
\end{example}

\begin{proposition}\label{proposition:2}
{\rm
Let $X=\{x, y\}$, $E=\{e_1,e_2\}$.
If $(X, {\cal T}, E)$ is a soft $T_1$-\Rmnum{2} space, then $(F,E)$ and $(G,E)$ belong to $\cal T$, where
$$F({e_1}) = F({e_2})=\{x\},\ \ G({e_1}) =G({e_2})= \{ y \}.$$
}
\end{proposition}
\noindent {\bf Proof.}
Since $(X, {\cal T}, E)$ is a soft $T_1$-\Rmnum{2} space and $P_{e_1}^ x\neq P_{e_2}^ x$, there exists a soft open set
$(H_1,E)$ such that $x\in H_1(e_1),x\notin H_1(e_2)$.
For $P_{e_1}^ x\neq P_{e_2}^ y$, there exists a soft open set
$(H_2,E)$ such that $x\in H_2(e_1),y\notin H_2(e_2)$.
Similarly, for $P_{e_1}^ x\neq P_{e_1}^ y$, there exists a soft open set
$(H_3,E)$ such that $x\in H_3(e_1),y\notin H_3(e_1)$.
Let $(H_4,E)=(H_1,E)\cap (H_2,E)\cap (H_3,E)$, we have $(H_4,E)\in {\cal T}$, where $H_4$ satisfies
$x\in H_4(e_1),y\notin H_4(e_1), x\notin H_4(e_2),y\notin H_4(e_2)$.
Thus $H_4(e_1)=\{x\}, H_4(e_2)=\emptyset$.

For $P_{e_1}^ x\neq P_{e_2}^ x$, there exists a soft open set $(H_5,E)$
such that $x\in H_5(e_2),x\notin H_5(e_1)$.
For $P_{e_1}^y\neq P_{e_2}^ x$, there exists a soft open set $(H_6,E)$
such that $x\in H_6(e_2),y\notin H_6(e_1)$.
Similarly, for $P_{e_2}^x\neq P_{e_2}^ y$, there exists a soft open set $(H_7,E)$
such that $x\in H_7(e_2),y\notin H_7(e_2)$. Let $(H_8,E)=(H_5,E)\cap (H_6,E)\cap (H_7,E)$, we have $(H_8,E)\in {\cal T}$ and
$H_8(e_1)=\emptyset,H_8(e_2)=\{x\}$.

Since $\cal T$ is closed on the operation union, we have $(H_4,E)\cup (H_8,E)=(F,E)\in {\cal T}$.
The same method can be used to prove the conclusion $(G,E)\in{\cal T}$.
\qed

\begin{remark}\label{remark:3}
{\rm
Example \ref{example:7} shows that a soft $T_1$-\Rmnum{1} space may not to be a soft $T_1$-\Rmnum{2} space.
Meanwhile it indicates that the converse of Proposition \ref{proposition:2} may not be true.
}
\end{remark}

By Proposition \ref{proposition:2} and Definition \ref{definition:13}, we can obtain the following conclusion immediately.
\begin{proposition}\label{proposition:7}
{\rm
If the initial universe set $X$ and the parameter set $E$ are binary sets,
then a soft $T_1$-\Rmnum{2} space $(X, {\cal T}, E)$ over $X$ is a soft $T_1$-\Rmnum{1} space.
}
\end{proposition}

\subsection{Soft $T_2$-\Rmnum{1} space and soft $T_1$-\Rmnum{2} space}
\begin{definition}\label{definition:15}{\rm (\cite{topological})}
{\rm
Let $(X, {\cal T}, E)$ be a soft topological space over $X$ and $x, y\in X$ such that $x\neq y$.
If there exist soft open sets $(F , E)$ and $(G, E)$ such that
$x\widehat{\in} (F , E)$ ,
$y\widehat{\in} (G, E)$ and $(F,E)\cap (G, E)=\widetilde{\emptyset}$,
then $(X, {\cal T}, E)$ is called a soft $T_2$-\Rmnum{1} space.
}
\end{definition}

\begin{definition}\label{definition:16}{\rm (\cite{Separation})}
{\rm
Let $(X, {\cal T}, E)$ be a soft topological space over $X$  and $P_a^x, P_b^y\in S\!P(S\!S_E(X))$ such that $P_a^x\neq P_b^y$.
If there exist soft open sets $(F , E)$ and $(G, E)$ such that
$P_a^x\widetilde{\in} (F , E)$ ,
$P_b^y\widetilde{\in} (G, E)$ and $(F,E)\cap (G, E)=\widetilde{\emptyset}$,
then $(X, {\cal T}, E)$ is called a soft $T_2$-\Rmnum{2} space.
}
\end{definition}

By definitions defined above, every soft $T_2$-\Rmnum{1} space is a soft $T_1$-\Rmnum{1} space,
and every soft $T_1$-\Rmnum{1} space is a soft $T_0$-\Rmnum{1} space.
Every soft $T_2$-\Rmnum{2} space is a soft $T_1$-\Rmnum{2} space,
and every soft $T_1$-\Rmnum{2} space is a soft $T_0$-\Rmnum{2} space(see \cite{topological,Separation}).

\begin{example}\label{example:8}
{\rm As shown in Example \ref{example:7}, we known that $(X, {\cal T}, E)$ is not a soft $T_2$-\Rmnum{2} topological space.

For $x$ and $y$, there exist two soft sets $(F,E)$ and $(G,E)$ such that
$x\widehat{\in} (F , E)$ ,
$y\widehat{\in} (G, E)$ and $(F,E)\cap (G, E)=\widetilde{\emptyset}$.
Then $(X, {\cal T}, E)$ is a soft $T_2$-\Rmnum{1} space.
This example shows that a soft $T_2$-\Rmnum{1} space may not to be a soft $T_2$-\Rmnum{2} space.
}
\end{example}

\begin{proposition}\label{proposition:3}
{\rm
Let $X=\{x, y\}$, $E=\{e_1,e_2\}$.
If $(X, {\cal T}, E)$ is a soft $T_1$-\Rmnum{2} space over $X$, then it is a soft $T_2$-\Rmnum{1} space.
}
\end{proposition}
\noindent {\bf Proof.}
By Proposition \ref{proposition:2}, we have $(F,E)$ and $(G,E)$ defined in Proposition \ref{proposition:2} belong to $\cal T$.
Thus $x\widehat\in (F, E), y\widehat\in(G, E)$ and $(F, E)\cap (G, E) = \widetilde\emptyset$.
Then $(X, {\cal T}, E)$ is a soft $T_2$-\Rmnum{1} space.
\qed

By Proposition \ref{proposition:3}, we have the following conclusion.
\begin{proposition}\label{proposition:8}
{\rm
If the initial universe set $X$ and the parameter set $E$ are binary sets,
then a soft $T_2$-\Rmnum{2} space $(X, {\cal T}, E)$ over $X$ is a soft $T_2$-\Rmnum{1} space.

}
\end{proposition}

\begin{remark}\label{remark:6}
{\rm
Proposition \ref{proposition:1}, Proposition \ref{proposition:7} and Proposition \ref{proposition:8}
only proved the case of binary sets.
In general, if $X$ and $E$ are finite sets or arbitrary sets, is this conclusion correct that
a soft $T_2$-\Rmnum{2}($T_1$-\Rmnum{2} and $T_0$-\Rmnum{2} respectively) space over $X$
is a soft $T_2$-\Rmnum{1} ($T_1$-\Rmnum{1} and $T_0$-\Rmnum{1} respectively) space?
}
\end{remark}

\subsection{Soft $T_3$-\Rmnum{1} space and soft $T_3$-\Rmnum{2} space}
\begin{definition}\label{definition:17}{\rm (\cite{topological})}
{\rm
Let $(X, {\cal T}, E)$ be a soft topological space over $X$,
$(G, E)$ be a soft closed set in $X$ and $x\in X$ such that $x\widehat\notin (G, E)$.
If there exist soft open sets $(F_1, E)$ and $(F_2, E)$ such that
$x\widehat\in (F_1, E), (G, E)\subseteq(F_2, E)$ and $(F_1, E)\cap (F_2, E) = \widetilde\emptyset$,
then $(X, {\cal T} , E)$ is called a soft type-\Rmnum{1} regular space.

If a soft type-\Rmnum{1} regular space $(X, {\cal T}, E)$ is a soft $T_1$-\Rmnum{1} space,
then it is called a soft $T_3$-\Rmnum{1} space.
}
\end{definition}

\begin{definition}\label{definition:18}{\rm (\cite{Separation})}
{\rm
Let $(X, {\cal T}, E)$ be a soft topological space over $X$,
$(G, E)$ be a nontrivial soft closed set in $X$ and $P_a^x\in S\!P(S\!S_E(X))$ such that $P_a^x\widetilde\notin (G, E)$.
If there exist soft open sets $(F_1, E)$ and $(F_2, E)$ such that
$P_a^x\widetilde\in (F_1, E), (G, E)\subseteq(F_2, E)$ and $(F_1, E)\cap (F_2, E) = \widetilde\emptyset$,
then $(X, {\cal T} , E)$ is called a soft type-\Rmnum{2} regular space.

If a soft type-\Rmnum{2} regular space $(X, {\cal T}, E)$ is a soft $T_1$-\Rmnum{2} space,
then it is called a soft $T_3$-\Rmnum{2} space.
}
\end{definition}

\begin{example}\label{example:9}
{\rm
As shown in Example \ref{example:6}, every soft open set is also a soft closed set in
the soft topological space $(X, {\cal T}, E)$.

For $x$ and a soft closed set $P_{e_1}^x$, we have $x\widehat\notin P_{e_1}^x$.
There is only one soft set ${\widetilde X}$ such that $x\widehat\in{\widetilde X}$.
While $P_{e_1}^x\cap {\widetilde X}\neq \widetilde\emptyset$.
Then $(X, {\cal T}, E)$ is not a soft type-\Rmnum{1} regular topological space.

For a soft point $P_{e_i}^x$, there is only a nontrivial soft closed set $ P_{e_j}^x(i\neq j)$ such that
$P_{e_i}^x\widetilde\notin P_{e_j}^x$.
We have $P_{e_i}^x\widetilde\in P_{e_i}^x, P_{e_j}^x\widetilde\in P_{e_j}^x$ and
$P_{e_i}^x\cap P_{e_j}^x=\widetilde\emptyset$.
Then $(X, {\cal T}, E)$ is a soft type-\Rmnum{2} regular space.

This example shows that a soft type-\Rmnum{2} regular space may not to be a soft type-\Rmnum{1} regular space.
}
\end{example}

\begin{example}\label{example:12}
{\rm
Let $X=\{x, y\}$, $E=\{e_1,e_2\}$ and ${\cal T}=\{ \widetilde\emptyset ,\widetilde X, P_{{e_1}}^x, (F,E)\}$, where
\begin{center}
$F(e_1)=\{y\}, F(e_2)=X$.
\end{center}
Then $(X, {\cal T}, E)$ is a soft topological space
and $\cal T$ is also the family of the soft closed sets in the soft topological space.

For $P_{{e_1}}^x\in {\cal T}$, we have $x\widehat\notin P_{{e_1}}^x$.
There exists only one soft set $\widetilde X$ such that $x\widehat\in \widetilde X$.
So there cannot be a soft open set $(G,E)$ which contains $P_{{e_1}}^x$ and satisfies
$(G,E)\cap \widetilde X=\widetilde\emptyset$.
Thus $(X, {\cal T}, E)$ is not a soft type-\Rmnum{1} regular space.

We know that $P_{{e_1}}^x$ and $(F,E)$ are two nontrivial soft closed sets in the soft topological space, and
$S\!P(S\!S_E(X))=\{P_{e_1}^{x}, P_{e_1}^{y}, P_{e_2}^{x}, P_{e_2}^{y}\}$.
There are four cases where the relationship $\widetilde\in$ between a soft point and
a nontrivial soft closed set does not hold:
$P_{{e_1}}^x\widetilde\notin(F,E)$, $P_{{e_2}}^x\widetilde\notin P_{e_1}^{x}$, $P_{{e_1}}^y\widetilde\notin P_{e_1}^{x}$ and
$P_{{e_2}}^y\widetilde\notin P_{e_1}^{x}$.
For example, for $P_{{e_1}}^y\widetilde\notin P_{e_1}^{x}$,
we have $P_{{e_1}}^y\subseteq (F,E)$, $P_{e_1}^{x}\widetilde\in P_{e_1}^{x}$
and $(F,E)\cap P_{e_1}^{x}=\widetilde\emptyset$. The other three cases can be shown similarly.
Then $(X, {\cal T}, E)$ is a soft type-\Rmnum{2} regular space.

Thus, if the sets $X$ and $E$ are binary sets,
a soft type-\Rmnum{2} regular space over $X$ may not to be a soft type-\Rmnum{1} regular space.
It is different from Proposition \ref{proposition:1}, Proposition \ref{proposition:7} and Proposition \ref{proposition:8}.
}
\end{example}

\begin{proposition}\label{proposition:4}
{\rm
Let $(X, {\cal T}, E)$ is a soft topological space over $X$.
If it is a soft type-\Rmnum{1} regular space, then it is a soft type-\Rmnum{2} regular space.
}
\end{proposition}
\noindent {\bf Proof.}
Let $P_a^x\in S\!P(S\!S_E(X))$, $(K,E)$ be a soft closed set and $P_a^x\widetilde\notin (K,E)$.
Then $x\widehat\notin (K,E)$. Since $(X, {\cal T}, E)$ is a soft type-\Rmnum{1} regular space,
there exist two soft open sets $(F,E)$ and $(G,E)$ such that
$x\widehat\in (F, E), (K, E)\subseteq(G, E)$ and $(F, E)\cap (G, E) = \widetilde\emptyset$.
Clearly that $P_a^x\widetilde\in (F, E)$. By Definition \ref{definition:18},
$(X, {\cal T}, E)$ is a soft type-\Rmnum{2} space.
\qed

\begin{example}\label{example:10}
{\rm Let $X=\{x, y, z\}$, $E=\{e_1,e_2\}$ and
${\cal T}=\{\widetilde\emptyset ,\widetilde X,(F_1,E ),(F_2,E ),(F_3,E),\\
(F_4, E), (F_5, E), (F_6, E)\}$,
where\\
${F_1}({e_1}) = \{x,y \},{F_1}({e_2})=\{x, y\}$,\\
${F_2}({e_1}) = \{ x,z \},{F_2}({e_2}) = \{ x,z \}$,\\
${F_3}({e_1}) = \{ y,z \},{F_3}({e_2}) = \{ y,z \}$,\\
${F_4}({e_1}) = \{ x \},{F_4}({e_2}) =\{x \}$,\\
${F_5}({e_1}) = \{ y \},{F_5}({e_2}) = \{ y \}$,\\
${F_6}({e_1}) = \{ z \},{F_6}({e_2}) = \{ z \}.$\\
Then $(X, {\cal T}, E)$ is a soft topological space.

For $x$ and $y$, there exist two soft open sets $(F_4,E)$ and $(F_5,E)$ such that
$x\widehat\in (F_4,E), \\ y\widehat\notin (F_4,E)$ and $y\widehat\in (F_5,E), x\widehat\notin (F_5,E)$.
For $x$ and $z$, or $y$ and $z$, there are also corresponding soft open sets to separate single point soft sets
generated by them.
Thus $(X, {\cal T}, E)$ is a soft $T_1$-\Rmnum{1} space.

For $x\in X$, there are three soft closed sets $(F_1,E)^{'}, (F_2,E)^{'}$ and $(F_4,E)^{'}$ that do not contain $x$
with respect to the relationship $\widehat\in$.
We explain that they meet the conditions of soft type-\Rmnum{1} regular spaces one by one.
For $x$ and $(F_1,E)^{'}$, there exist two soft open sets $(F_4,E)$ and $(F_6,E)$ such that
$x\widehat\in(F_4,E), (F_1,E)^{'}\subseteq (F_6,E)$ and  $(F_4,E)\cap (F_6,E)=\widetilde\emptyset$.
For $x$ and $(F_2,E)^{'}$, there exist two soft open sets $(F_4,E)$ and $(F_5,E)$ such that
$x\widehat\in(F_4,E), (F_2,E)^{'}\subseteq (F_5,E)$ and  $(F_4,E)\cap (F_5,E)=\widetilde\emptyset$.
For $x$ and $(F_4,E)^{'}$, there exist two soft open sets $(F_4,E)$ and $(F_3,E)$ such that
$x\widehat\in(F_4,E), (F_4,E)^{'}\subseteq (F_3,E)$ and  $(F_4,E)\cap (F_3,E)=\widetilde\emptyset$.

Similarly, for $y,z\in X$ and the corresponding soft closed sets which do not contain $y$ or $z$,
they also can be separated by some soft open sets.
Thus $(X, {\cal T}, E)$ is a soft type-\Rmnum{1} regular space.
And then $(X, {\cal T}, E)$ is a soft $T_3$-\Rmnum{1} space.

For $P_{e_1}^x$ and  $P_{e_2}^x$, we have $P_{e_1}^x\neq P_{e_2}^x$.
But every soft open set $(F,E)$ satisfying $x\in F(e_1)$ has $x\in F(e_2)$.
That is to say, if $P_{e_1}^x\widetilde\in (F, E)$, then $P_{e_2}^x\widetilde\in (F, E)$.
So $(X, {\cal T}, E)$ is not a soft $T_1$-\Rmnum{2} space.
And then $(X, {\cal T}, E)$ is not a soft $T_3$-\Rmnum{2} space.
}
\end{example}

\begin{remark}\label{remark:7}
{\rm
Example \ref{example:10} shows that a soft $T_3$-\Rmnum{1} space may not to be a soft $T_3$-\Rmnum{2} space.
}
\end{remark}

\subsection{Soft $T_4$-\Rmnum{1} space and soft $T_4$-\Rmnum{2} space}
\begin{definition}\label{definition:19}{\rm (\cite{topological})}
{\rm
Let $(X, {\cal T}, E)$ be a soft topological space over $X$,
$(F, E)$ and $(G, E)$ be a soft closed set over $X$ such that $(F, E)\cap (G, E)=\widetilde\emptyset$.
If there exist soft open sets $(F_1, E)$ and $(F_2, E)$ such that
$(F, E)\subseteq (F_1, E), (G, E)\subseteq(F_2, E)$ and $(F_1, E)\cap (F_2, E) = \widetilde\emptyset$,
then $(X, {\cal T} , E)$ is called a soft normal space.
}
\end{definition}

\begin{definition}\label{definition:20}{\rm (\cite{topological,Separation})}
{\rm Let $(X, {\cal T}, E)$ be a soft topological space over $X$.
Then it is said to be a soft $T_4$-\Rmnum{1} space if it is soft normal space and soft $T_1$-\Rmnum{1} space.
If a soft normal space $(X, {\cal T}, E)$ is a soft $T_1$-\Rmnum{2} space,
then it is called a soft $T_4$-\Rmnum{2} space.
}
\end{definition}

\begin{example}\label{example:11}
{\rm
As described in Example \ref{example:10}, the soft topological space $(X, {\cal T}, E)$
has six soft closed sets $(F_i,E)^{'}(i=1,2,\cdots,6)$.
They form six pairs of disjoint soft closed sets: $(F_1,E)^{'}$ and $(F_2,E)^{'}$, $(F_1,E)^{'}$ and $(F_3,E)^{'}$,
$(F_1,E)^{'}$ and $(F_6,E)^{'}$, $(F_2,E)^{'}$ and $(F_3,E)^{'}$, $(F_2,E)^{'}$ and $(F_5,E)^{'}$,
$(F_3,E)^{'}$ and $(F_4,E)^{'}$.

For the pair of soft closed set $(F_1,E)^{'}$ and $(F_2,E)^{'}$, there exist soft open sets $(F_5,E)$ and $(F_6,E)$
such that $(F_1,E)^{'}\subseteq (F_6,E)$, $(F_2,E)^{'}\subseteq (F_5,E)$ and $(F_5, E)\cap (F_6, E) = \widetilde\emptyset$.
For the pair of soft closed sets $(F_1,E)^{'}$ and $(F_6,E)^{'}$, there exist soft open sets $(F_1,E)$ and $(F_6,E)$
such that $(F_1,E)^{'}\subseteq (F_6,E)$, $(F_6,E)^{'}\subseteq (F_1,E)$ and $(F_1, E)\cap (F_6, E) = \widetilde\emptyset$.
Similarly, pairs of other disjoint soft closed sets also can be proved. Thus $(X, {\cal T} , E)$ is a soft normal space.

Since
$(X, {\cal T} , E)$ is not a soft $T_1$-\Rmnum{2} space, it is not a soft $T_4$-\Rmnum{2} space.
But it is a soft $T_1$-\Rmnum{1} space, and then it is a soft $T_4$-\Rmnum{1} space.

This example shows that a soft $T_4$-\Rmnum{1} space may not to be a soft $T_4$-\Rmnum{2} space.
}
\end{example}

\begin{remark}\label{remark:8}
{\rm
If the initial universe set $X$ and the parameter set $E$ are binary sets,
a soft $T_4$-\Rmnum{2} space over $X$ is a soft $T_4$-\Rmnum{1} space by discussion above.
}
\end{remark}

\section{Conclusion}\label{sec:con}
In this paper, we gave some examples to verify differences between two types of separation axioms in soft topological
spaces. We have known that a soft $T_0$-\Rmnum{1}($T_1$-\Rmnum{1}, $T_2$-\Rmnum{1}, $T_3$-\Rmnum{1} and $T_4$-\Rmnum{1}
respectively) space over $X$
may not to be a soft $T_0$-\Rmnum{2}($T_1$-\Rmnum{2}, $T_2$-\Rmnum{2}, $T_3$-\Rmnum{2} and $T_4$-\Rmnum{2} respectively) space.
If the initial universe set $X$ and the parameter set $E$ are binary sets, the converse are true except for the property $T_3$.
We also showed that a soft type-\Rmnum{1} regular space is a soft type-\Rmnum{2} regular space.

\section*{Acknowledgments}
This work is supported by National Science Foundation of China (Grant No.12071188).


\begin{thebibliography}{99} \addtolength{\itemsep}{-1.5ex}

\bibitem{first} D. Molodtsov, Soft set theory-First results, \textit{Comput. Math. Appl. }\textbf{37}(1999), 19-31.

\bibitem{soft} P.K. Maji, R. Biswas, A.R. Roy, Soft set theory, \textit{Comput. Math. Appl.}\textbf{45}(2003), 555-562.

\bibitem{operation} M. Irfan Ali, F. Feng,  X.Y. Liu, W.K. Min, M. Shabir,
On some new operations in soft set theory, \textit{Comput. Math. Appl.}\textbf{57} (2009), 1547--1553.

\bibitem{category}
R.A. Borzooei, M. Mobini, M.M, Ebrahimi, The category of soft sets,
\textit{Journal of Intelligent \& Fuzzy Systems} \textbf{28}(2015),157-167.

\bibitem{the group} H. Akta\c{s}, N. \c{C}a\u{g}man, Soft sets and soft groups,
\textit{Inform. Sci.}\textbf{177}(2007), 2726-2735.

\bibitem{algebra} Y.B. Jun, K.J. Lee, J.M. Zhan, Soft p-ideal of soft BCK/BCI-algebras,
\textit{Comput. Math. Appl.} \textbf{58}(2009), 2060-2068.

\bibitem{ideal} Y.B. Jun, C.H. Park, Applications of soft sets in ideal theory of BCK/BCI-algebras,
\textit{Inform. Sci.} \textbf{178}(2008), 2466-2475.

\bibitem{Aktas} H. Akta\c{s}, Some algebraic applications of soft sets,
\textit{Applied Soft Computing }\textbf{28}(2015), 327--331.

\bibitem{semirings} F. Feng, Y.B. Jun, X.Z. Zhao, Soft semirings, \textit{Comput. Math. Appl.}
\textbf{56}(2008), 2621-2628.

\bibitem{rough} F. Feng, X.Y. Liu, V.Leoreanu-Fotea, Y.B. Jun, Soft sets and soft rough sets,
\textit{Inform. Sci.} \textbf{181}(2011), 1125-1137.

\bibitem{topology}
N. \c{C}a\u{g}man, S. Karatas, S. Enginoglu, Soft topology,
\textit{Comput. Math. Appl.}\textbf{62}(2011), 351--358.

\bibitem{topological} M.Shabir, M.Naz , On soft topological spaces,
\textit{Computers and Mathematics with Applications} \textbf{61}(2011), 1786-1799.


\bibitem{Sabir} Sabir Hussain, Bashir Ahmad, Soft separation axioms in soft topological spaces,
\textit{Hacettepe Journal of Mathematics and Statistics} \textbf{44}(2015), 559--568.

\bibitem{note} W.K. Min, A note on soft topological spaces,
 \textit{Computers and Mathematics with Applications} \textbf{62}(2011), 3524-3528.

\bibitem{remark} I. Zorlutuna, M. Akdag, W.K. Min, Remarks on soft topological spaces,
\textit{Annals of fuzzy mathematics and informatics} \textbf{3}(2012), 171-185.

\bibitem{Arif1}
Arif Mehmood Khattak, Gulzar Ali Khan, Muhammad Ishfaq, Fahad Jamal,
Characterization of soft $\alpha$-separation axioms and soft $\beta$-separation axioms in
soft single point spaces and in soft ordinary spaces,
\textit{Journal of New Theory} \textbf{19}(2017), 63--81.

\bibitem{Arif2}
A.M. Khattak, A. Zaighum, G.A. Khan, I. Ahmed, S. Khan,
Soft separation axioms in soft single point spaces and in soft ordinary spaces,
\textit{Asian Journal of Mathematics and Computer Research} \textbf{23}(2018), 175--192.

\bibitem{Orhan} Orhan G\"{o}\c{c}\"{u}r, Abdullah Kopuzlu, Some new properties on soft separation axioms,
\textit{Annals of Fuzzy Mathematics and Informatics} \textbf{9}(2015), 421--429.


\bibitem{Sabir2} Sabir Hussain, Bashir Ahmad, On some structures of soft topology,
\textit{Mathematical Sciences} \textbf{6}(2012), 7 pages.

\bibitem{Separation} J.L. He, Separation properties in soft topological spaces,
\textit{Pure and Applied Mathematics}(in Chinese) \textbf{33}(2017), 1419--151.

\bibitem{zhang} X.W. Zhang Xiongwei, Further study on soft separation axioms,
\textit{Computer Engineering and Applications} \textbf{49}(2013), 48--50.


\bibitem{fuzzy} P.K. Maji, R. Biswas, A.R. Roy, Fuzzy Soft Sets,
\textit{Journal of Fuzzy Mathematics} \textbf{9}(2001), 589--602.

\bibitem{continuous}
X.C. Guan, Y.M. Li, F. Feng, A new order relation on fuzzy soft sets
and its application, \textit{Soft Computing} \textbf{17}(2013), 63--70.

\bibitem{interval}
M. Son, Interval-valued Fuzzy Soft Sets,
\textit{Journal of Korean Institute of Intelligent Systems }\textbf{4}(2007), 557--562.


\bibitem{reduction} D. Chen, E.C.C. Tsang, D.S. Yeung, X. Wang, The parametrization reduction of soft sets
and its applications, \textit{Comput. Math. Appl.} \textbf{49}(2005), 757-763.

\bibitem{Deli} I. Deli, N. \c{C}a\u{g}man,
Intuitionistic fuzzy parameterized soft set theory and its decision making,
\textit{Applied Soft Computing }\textbf{28}(2015), 109-113.

\bibitem{banghe}B.H. Han, S.L. Geng, Pruning method for optimal solutions of
int$^m$-int$^n$ decision making scheme \textit{ Eur. J. Oper. Res.} \textbf{231} (2013), 779-783.

\bibitem{making} H. Akta\c{s}, N. \c{C}a\u{g}man,
Soft decision making methods based on fuzzy sets and soft sets,
\textit{Journal of Intelligent \& Fuzzy Systems} \textbf{30}(5)(2016), 2797-2803.


\end{thebibliography}
\end{document}